\documentclass[11pt]{article}
\usepackage{epsfig}
\usepackage{t1enc}
    \usepackage[latin1]{inputenc}
    \usepackage[english]{babel}
        \usepackage{latexsym}
\usepackage{amssymb}
\usepackage{euscript}

\newcommand{\bs}{\textbf{S}}
\newcommand{\bo}{\textbf{0}}
\newcommand{\bb}{\textbf{B}}

\def\begg{\begin{equation}}
\def\endd{\end{equation}}

\newcommand{\bp}{{\bf P}}
\newcommand{\ep}{{\epsilon}}
\begin{document}
\setlength{\arraycolsep}{.136889em}
\renewcommand{\theequation}{\thesection.\arabic{equation}}
\newtheorem{thm}{Theorem}[section]
\newtheorem{propo}{Proposition}[section]
\newtheorem{lemma}{Lemma}[section]
\newtheorem{corollary}{Corollary}[section]
\newtheorem{remark}{Remark}[section]

\smallskip

\centerline{\Large\bf SOME LIMIT THEOREMS FOR HEIGHTS OF RANDOM WALKS}
\smallskip
\centerline{\Large\bf ON A SPIDER}

\bigskip\bigskip
\bigskip\bigskip
\centerline{\it Dedicated to  the memory of Marc Yor}

\bigskip\bigskip
\bigskip\bigskip

\bigskip\bigskip
\noindent
\textbf{Endre Cs\'{a}ki}
\newline
Corresponding author. Alfr\'ed R\'enyi Institute of Mathematics, Hungarian
Academy of Sciences, Budapest, P.O.B. 127, H-1364, Hungary. E-mail address:
csaki.endre@renyi.mta.hu
\newline
Phone: 361-483-8306,\, \, \, Fax: 361-483-8333

\bigskip
\noindent
\textbf{Mikl\'os Cs\"org\H{o}}
\newline
School of Mathematics and Statistics, Carleton University, 1125
Colonel By Drive, Ottawa, Ontario, Canada K1S 5B6. E-mail address:
mcsorgo@math.carleton.ca

\bigskip
\noindent
\textbf{Ant\'{o}nia F\"{o}ldes}
\newline
Department of Mathematics, College of Staten
Island, CUNY, 2800 Victory Blvd., Staten Island, New York 10314,
U.S.A.  E-mail address: antonia.foldes@csi.cuny.edu

\bigskip
\noindent
\noindent
\textbf{P\'al R\'ev\'esz}
\newline
\noindent
Institut f\"ur
Statistik und Wahrscheinlichkeitstheorie, Technische Universit\"at
Wien, Wiedner Hauptstrasse 8-10/107 A-1040 Vienna, Austria.
E-mail address: revesz.paul@renyi.mta.hu
\bigskip

\medskip
\noindent{\bf Abstract}\newline
A simple symmetric random walk is considered on a spider that is a collection
of half lines (we call them legs) joined at the origin. We establish a strong
approximation of this random walk by the so-called Brownian spider.
Transition probabilities are studied, and for a fixed number of legs we
investigate how high the walker and the Brownian motion can go on the legs
in $n$ steps. The heights on the legs are also investigated when the number
of legs goes to infinity.


\medskip

\vspace{.1cm}

\medskip

\section{Introduction }
\renewcommand{\thesection}{\arabic{section}} \setcounter{equation}{0}
\setcounter{thm}{0} \setcounter{lemma}{0}

Paraphrasing Harrison and  Shepp  \cite{HS}, in 1965 It\^{o} and McKean
(\cite{IM}, Section 4.2, Problem 1) introduced a simple but intriguing
diffusion process that they called skew Brownian motion, that was revisited
by Walsh \cite{WAL} in 1978. Walsh  introduced it as a Brownian motion with
excursions around zero in random directions on the plane. The random
directions are values of a random variable in $[0,2\pi)$ that are independent
for different excursions with a constant value during each excursion. This
"definition" can be made precise as, e.g., in Barlow, Pitman and Yor
\cite{BPY89}. This motion is now called Walsh's Brownian motion.

Following Barlow et al. \cite{BPY2} and Example 1 in Evans and Sowers \cite{ES},
we consider a version of Walsh's Brownian motion which lives on $N$ semi-axes
on the plane, called legs from now on, that are joined at the origin, the
so-called Brownian spider, or Walsh's spider. Loosely speaking, this motion
performs a regular Brownian motion on each one of the legs and, when it
arrives to the origin, it continues its motion on any of the $N$ legs
with a given probability. Thus, one can construct the Brownian spider by
independently putting the excursions from zero of a standard Brownian motion
on the $j$-th leg of the spider with probability $p_j$, $j=1,2,\ldots,N$
with $\sum_{j=1}^Np_j=1$. For a formal definition of the Brownian spider
along these lines we refer to Section 2, (\ref{sdef}). In the special case of
$p_j=1/N,\, j=1,2,\ldots,N$, Papanicolaou et al. \cite{PPL} studied the exit
time of this motion from specific sets and introduced a generalized arc-sine
law as well, concerning the time spent globally on the legs. This question
was further investigated in the elegant paper of  Vakeroudis  and Yor \cite{VY}.

A natural discrete counterpart of this motion is that of random walks on a
spider, i.e., replacing the Brownian motions with simple symmetric random
walks on the legs. Hajri \cite{Ha} studied such discrete versions as
approximations of the Brownian spider, proving their weak convergence to the
latter in a more general context of discrete approximations that are related
to Walsh's Brownian motion. As to the weak convergence in hand, he showed
that it can be deduced from the special case of $N=2$ converging to a skew
Brownian motion. Harrison and Shepp \cite{HS} reviewed the construction of
a skew Brownian motion from its scale and speed measure and considered it to
be a solution of a particular stochastic equation. Completing a random walk
result of \cite{HS}, Cherny et al. \cite{CSY} concluded weak convergence of
skew random walk to skew Brownian motion. For further discussions and
references we refer to Lejay \cite{Le}.

For the sake of studying random walks on the just mentioned spider,
we proceed with concrete definitions in this regard.

Put $\textbf{SP}(N)=(V_N,E_N),$ where, with $i=\sqrt{-1}$,
\begg V_N=\left\{v_N(r,j)=r\exp\left(\frac{2\pi ij}{N}\right),
\quad r=0,1,...,\quad j=1,...,N\right\} \label{veen}\endd
\noindent
is the set of vertices of $\textbf{SP}(N)$, and
$$E_N=\{e_N(r,j)=(v_N(r,j),v_N(r+1,j)),  \quad \,\, r=0,1,...,\,\,\quad
j=1,...,N \},$$
\noindent
is the set of edges of $\textbf{SP}(N)$. We will call the graph
$\textbf{SP}(N)$ a spider with $N$ legs. The vertex
$$v_N(0):=v_N(0,1)=v_N(0,2)=...=v_N(0,N)$$
\noindent
is called the body of the spider, and
$\{v_N(1,j), v_N(2,j),...\}$ is the $j$-th leg of the spider.
When the number of legs $N$ is fixed, we will suppress it in the notation
and, instead of $v_N(r,j)$ or $v_N(0)$, we will simply write $v(r,j)$
or $v(0)=\bo$, whenever convenient.

In this paper we consider a random walk $\bs_n,\,\, n=1,2 \ldots$, on
$\textbf{SP}(N)$ that starts from the body of the spider, i.e.,
$\bs_0=v_N(0)=\bo,$ with the following transition probabilities:
$${\bf P} (\bs_{n+1}=v_N(1,j)|\bs_n=v_N(0))=p_j,\quad j=1,...,N, $$
with
$$\sum_{j=1}^N p_j=1,$$
and, for $r=1,...,\quad j=1,...,N$,
$${\bf P} (\bs_{n+1}=v_N(r+1,j)|\bs_n=v_N(r,j))={\bf P}
(\bs_{n+1}=v_N(r-1,j)|\bs_n=v_N(r,j))=\frac{1}{2}.$$

The random walk $\bs_n$ on spider $\textbf{SP}(N)$ can be constructed from a
simple symmetric random walk $S(n)$, $n=0,1,\ldots$ on the line as follows.
Consider the absolute value $|S(n)|,\, n=1,2,\ldots$, that consists of
infinitely many excursions from zero, denoted by $G_1,G_2,\ldots$.
Put these excursions, independently of each other, on leg $j$ of the spider
with probability $p_j$, $j=1,2,\ldots,N$. Thus we obtain the first $n$ steps
of the spider walk $\bf{S}_\cdot$, as above, from the first $n$ steps of the
random walk $S(\cdot)$.

We denote the Brownian spider on $\textbf{SP}(N)$, as described in the
second paragraph above, by ${\bb}(t),\, t\geq 0$,  that also starts from
the body of the spider, i.e., ${\bb}(0)=v_N(0)$.

In his book R\'ev\'esz \cite{R13} discussed the spider walk above in the case
when $p_j=1/N$, and the number of legs of the spider goes to infinity. In
our just introduced definitions, we followed the latter book but allow the
walker to select the legs with possibly unequal probabilities. In particular,
one can construct this spider walk by independently putting the excursions
from zero of a simple symmetric random walk on the $j$-th leg of the spider
with probability $p_j$ as above. Hence, in what follows, we will frequently
make use of arguments in terms of the usual simple symmetric random walk on
the line. In view of this, in the sequel, $\bs_n$ will stand for spider
walk, and $S(n)$ for a simple symmetric random walk on the line with respective
probabilities denoted by ${\bf P}$ and $P$.

In our Section 2 we establish a strong invariance principle for
approximating the spider walk ${\bs}_n$ by the Brownian spider
${\bb}(n)$, keeping $N$ fixed. In Section 3 we investigate the transition
probabilities, while in Section 4 we discuss  how high the random walk can go
on a spider with $N$ legs, where  $N$ is still fixed. The last section,
Section 5, is devoted to studying the probability that the walk goes up to
certain heights simultaneously on all legs when the number of legs are
increasing.

\section{Strong approximations}
\renewcommand{\thesection}{\arabic{section}} \setcounter{equation}{0}
\setcounter{thm}{0} \setcounter{lemma}{0}

The Brownian spider can be constructed from a standard Brownian motion
$\{B(t),\, t\geq 0\}$ on the line as follows. The process
$\{|B(t)|,\, t\geq 0\}$ has a countable number of excursions from zero, and
let $J_1,J_2,\ldots$ denote a fixed enumeration of its excursion intervals
away from zero. Then, for any $t>0$ for which $B(t)\neq 0$, we have that
$t\in J_m$ for one of the values of $m=1,2,\ldots$.

Extend the definition of $v_N(r,j)$ given in (\ref{veen}) to all positive
values of $r$, i.e.,
$$
v_N(r,j)=r\exp\left(\frac{2\pi ij}{N}\right).
\quad r\geq 0,\quad j=1,...,N,$$
Thus, $v_N(r,j)$ is the $j$-th leg of the spider. Let
$\kappa_m, \, m=1,2,\ldots,$ be i.i.d. random variables, independent of $B$
with
$$
P(\kappa_m=j)=p_j,\quad j=1,2,\ldots, N.
$$
We now construct the Brownian spider $\{{\bb}(t),\, t\geq 0\}$ by putting the
excursion whose interval is $J_m$, to leg $\kappa_m$ on the spider
$\textbf{SP}(N)$. Hence we can define the Brownian spider as discussed in
paragraph 2 of our Introduction by
\begg
{\bb}(t):=\sum_{m=1}^\infty I\{t\in J_m\}v_N(|B(t)|,\kappa_m), \quad
{\rm if}\quad B(t)\neq 0, \label{sdef}
\endd
and
$$
{\bb}(t):=v_N(0)=\bo,\quad {\rm if}\quad B(t)=0,
$$
where $I\{...\}$ is the indicator function.

This definition of the Brownian spider $\{{\bb}(t),\, t\geq 0\}$ is an
analogue of that of a skew Brownian motion given in Appuhamillage et al.
\cite{AB11}. In this regard, we may also refer to Revuz and Yor, Exercise
2.16, Chap XII in \cite{RY}. We note in passing that the Brownian spider
with $N=2$ is equivalent to the skew Brownian motion.

Moreover, define the distance on $\textbf{SP}(N)$ by
$$
|v_N(x,j)-v_N(y,j)|=|x-y|, \quad j=1,\ldots,N
$$
$$
|v_N(x,j)-v_N(y,k)|=x+y, \quad j,k=1,\ldots,N,\, \, j\neq k.
$$

First we mention the weak convergence result of Hajri \cite{Ha}.
\begin{thm} Let ${\bf S}(t),\, t\geq 0,$ be the linear interpolation of
${\bf S}_n,\, n=0,1,\ldots$.
Then
$$
\left\{\frac{{\bf S}(nt)}{\sqrt{n}}, t\geq 0\right\}\to \{{\bf B}(t),\,
t\geq 0\}
$$
weakly on $C[0,\infty)$, as $n\to\infty$.
\end{thm}

Our strong approximation result, that also contains Theorem 2.1, reads as
follows.
\begin{thm} On a rich enough probability space one can define a Brownian
spider $\{{\bf B}(t),\, t\geq 0\}$ and a random walk
$\{{\bf S}_n,\, n=0,1,2,\ldots\}$, both on ${\bf SP}(N)$, and both
selecting their legs  with the same  probabilities
$p_j,\, j=1,2,\ldots,N$ so that, as $n\to\infty$, we have
$$
|{\bf S}_n-{\bf B}(n)|=O((n\log\log n)^{1/4}(\log n)^{1/2}) \quad a.s.
$$
\end{thm}
{\bf Proof}. Start with a Skorokhod embedding for $B(\cdot)$ and $S(\cdot),$
i.e., define
\begin{eqnarray*}
\tau_1&=&\inf\{t>0: \, |B(t)|=1\},\\
\tau_2&=&\inf\{t>\tau_1:\, |B(t)-B(\tau_1)|=1\},\\
&\ldots&\\
\tau_{i+1}&=&\inf\{t>\tau_i:\, |B(t)-B(\tau_i)|=1\}.
\end{eqnarray*}
Then $\{S(n):=B(\tau_n),\, n=1,2,\ldots\}$ is a simple symmetric random
walk on the line and, as $n\to\infty$, we have
$$
|S(n)-B(n)|=|B(\tau_n)-B(n)|=O((n\log\log n)^{1/4}(\log n)^{1/2})\quad a.s.
$$
The latter Skorokhod embedding of $B$ and $S$ is a special case of Theorem 1.5
of Strassen \cite{STR}, (cf. also R\'ev\'esz \cite{R13}, Theorem 6.1 when $S$
is a simple symmetric random walk).

Now construct ${\bb}(t)$ from $B(t)$ as described above. It is clear from
Skorokhod construction that an excursion of $S(\cdot)$ lies entirely within
its corresponding excursion of $B(\cdot)$, so construct ${\bs}_n$ by putting
this excursion on the same leg as the corresponding excursion of $B(\cdot)$.
We note in passing that small excursions of the underlying Brownian motion,
namely those that do not reach 1, are not needed for the construction of
$\{{\bs}_n:={\bb}(\tau_n),\, n=1,2,\ldots\}$, i.e., for that of a random walk
on $\textbf{SP}(N)$.

Consider now ${\bb}(n)$ and ${\bb}(\tau_n)$ when they are on the same leg.
Then, as $n\to\infty$,
$$
|{\bs}_n-{\bb}(n)|=|{\bb}(\tau_n)-{\bb}(n)|=
|B(\tau_n)-B(n)|=O((n\log\log n)^{1/4}(\log n)^{1/2})\quad a.s.
$$
However, when ${\bb}(n)$ and ${\bb}(\tau_n)$ are on different legs, then
$$
|{\bs}_n-{\bb}(n)|=|{\bb}(\tau_n)-{\bb}(n)|=|B(\tau_n)|+|B(n)|.
$$
But in this case there is a point $c_n$ between $n$ and $\tau_n$, where
$B(c_n)=0$, with
$|n-c_n|\leq |n-\tau_n|$ and $|\tau_n-c_n|\leq |n-\tau_n|$.
Since $\tau_1,\tau_2-\tau_1,\tau_3-\tau_2,\ldots$ is a sequence of random
variables with mean 1 and variance 1 (cf., e.g., page 54 in R\'ev\'esz
\cite{R13}), by the law of the iterated logarithm (LIL), as $n\to\infty$, we
get that
$$|\tau_n-n|=O(n\log\log n)^{1/2}=:a_n.$$
Now applying the Wiener large increments result of Cs\"org\H{o}-R\'ev\'esz
\cite{CR79}, (see also page 30 in \cite{CR81}), as $n\to\infty$, we obtain

\begin{eqnarray*}
|B(n)|&=&|B(n)-B(c_n)|\leq \sup_{\,0\leq s\leq a_n}
|B(n-s)-B(n)|+\sup_{\,0\leq s\leq a_n} |B(n+s)-B(n)|\\
&=& O((n\log\log n)^{1/4}(\log n)^{1/2})\quad a.s.,
\end{eqnarray*}
and similarly
$$
|B(\tau_n)|=|B(\tau_n)-B(c_n)|=O((n\log\log n)^{1/4}(\log n)^{1/2})\quad a.s.
$$
This completes the proof of the Theorem 2.2. $\Box$

\section{Transition probabilities}
\renewcommand{\thesection}{\arabic{section}} \setcounter{equation}{0}
\setcounter{thm}{0} \setcounter{lemma}{0}
We assume throughout that $\bs_0=\bo.$ Clearly, we have
$\displaystyle{{\bf P} (\bs_{2n}=\bo)=P(S(2n)=0)}.$

\begin{thm} For $i\geq 1,j \geq 1$ integers
\begin{eqnarray*}
&{\rm (i)}&{\bf P}({\bf S}_{2n+2k}=v(2j,\ell)|{\bf S}_{2k}={\bf 0})
=2p_{\ell}\,P(S(2n)=2j),
\quad j\leq n\\
&{\rm (ii)}&
{\bf P}({\bf S}_{2k+2n}=v(2i,\ell^*)|{\bf S}_{2k}=v(2j,\ell))
=2p_{\ell^*}P(S(2n)=2(j+i)),\quad i+j \leq n,\quad \ell\neq\ell^*\\
&{\rm (iii)}& {\bf P}({\bf S}_{2n+2k}=v(2i,\ell)|
{\bf S}_{2k}=v(2j,\ell))\\
&&=P(S(2n)=2(j-i))-(1-2p_{\ell})P(S(2n)=2(j+i)),
\quad  |i-j|\leq n.
\end{eqnarray*}
\end{thm}

\noindent\textbf{Proof}:
It is  well-known that  for the simple symmetric random walk we have

\begg
P(S(2n)=2k)=\frac{1}{2^{2n}} {2n\choose n+k},
\endd
and, for any integer $k\geq 1,$ we have from the ballot theorem, that
\begg
P(S(1)>0,\,S(2)>0,..., S(2n-1)>0,\,S(2n)=2k)=\frac{k}{n}
\frac{1}{2^{2n}}{2n \choose n+k}.    \label{bal}
\endd
\noindent
Partitioning according to the time of the last return to the origin,
using (\ref{bal}) and taking into account that the probability of the next
step after the last return to zero in our context is $p_\ell$ instead of
$1/2$, we get
\begin{eqnarray*}
{\bf P}(\bs_{2n}&=&v(2j,\ell))\\
=\sum_{m=0}^{n-1}{\bf P}(\bs_{2m}&=&0)2p_{\ell}P(S(1)>0,S(2)>0,
\ldots, S( 2(n-m)-1)>0,S(2(n-m))=2j )\\
=\sum_{m=0}^{n-1}{\bf P}(\bs_{2m}&=&0)p_{\ell}\,\frac{j}{n-m}\,
\frac{2}{2^{2(n-m)}}{2n-2m \choose n-m+j}\\
=2p_{\ell} \sum_{m=0}^{n-1}P(S(2m)&=&0)\,P(S(1)\neq 0,S(2)\neq 0,
\dots, S( 2(n-m)-1)\neq 0, S(2(n-m))=2j)\\
=2p_{\ell} P(S(2n)&=&2j),
\end{eqnarray*}
which proves (i).

For $\ell \neq \ell^*,$   partitioning again according to the last visit to
the origin, we arrive at
\begin{eqnarray*}
{\bf P}(\bs_{2k+2n}&=&v(2i,\ell^*)|\bs_{2k}=v(2j,\ell))\\
=\sum_{m=i}^{n-j}P(S(0)&=&2j, S(2n-2m)=0)\,2p_{\ell^*}\,
P(S(0)=0, S(1)>0,\ldots, S(2m-1)>0, S(2m)=2i)\\
&=&\frac{1}{2^{2n}}2p_{\ell^*}
\sum_{m=i}^{n-j}{2n-2m \choose n-m+j}{2m \choose m+i }\frac{i}{m}=
2p_{\ell^*} P(S(2n)=2(j+i)),
\end{eqnarray*}
which proves (ii).

Finally to prove (iii), observe that any path from $v(2j,\ell) $ to
$v(2i,\ell)$ either crosses the origin or not. For the transition with
crossing the origin we have
\begin{eqnarray*}
{\bf P}(\bs_{2k+2n}&=&v(2i,\ell)|\bs_{2k}=v(2j,\ell),
\bs_{2k+2m}=\bo\,\,{\rm for \,\,some}\,\, j \leq m \leq n-i)\\
={\bf P}(\bs_{2k+2n}&=&v(2i,\ell^*)|\bs_{2k}=v(2j,\ell))
\end{eqnarray*}
with the understanding that here leg $\ell^*$ is actually leg $\ell$,
and hence the probability gained in (ii) should be used with
$p_{\ell^*}=p_{\ell}.$ In the case when the transition happens without
crossing the origin, the corresponding probability can be calculated just
like for a simple symmetric walk, using the reflection principle. Thus

\begin{eqnarray*}
{\bf P}( \bs_{2k+2n}&=&v(2i,\ell)|\bs_{2n}=v(2j,\ell))\\
=2p_{\ell} P(S(2n)&=&2(j+i))+P(S(2n)=2(j-i))-P(S(2n)=2(j+i))\\
=P(S(2n)&=&2(j-i))-(1-2\,p_{\ell})\,P(S(2n)=2(j-i)).
\end {eqnarray*}
 proving (iii). $\Box$

Recall that by the local central limit theorem, as $n\to\infty$, we have
$$P(S(2n)=2k)\sim \frac{1}{\sqrt{\pi n}}e^{-\frac{k^2}{n}}$$
if $k/\sqrt{n}$ is bounded. Hence, via Theorem 3.1, we  obtain the
corresponding limit theorem for transition probabilities as follows.
\begin{thm}
\begin{eqnarray*}
&{\rm (i)}&
\lim_{n\to\infty}\sqrt{n}{\bf P}
({\bf S}_{2[nt]+2k}=v(2[y\sqrt{n}],\ell)|{\bf S}_{2k}={\bf 0})
=\frac{2p_{\ell}}{\sqrt{\pi t}}e^{-y^2/t},\\
&{\rm (ii)}&
\lim_{n\to\infty}\sqrt{n}{\bf P}({\bf S}_{2[nt]+2k}=
v(2[y\sqrt{n}],\ell^*)|{\bf S}_{2k}=v(2[x\sqrt{n}],\ell))
=\frac{2p_{\ell^*}}{\sqrt{\pi t}}e^{-(x+y)^2/t},\quad \ell\neq\ell^*,\\
&{\rm (iii)}& \lim_{n\to\infty}\sqrt{n}
{\bf P}({\bf S}_{2[nt]+2k}=v(2[y\sqrt{n}],\ell)|
{\bf S}_{2k}=v(2[x\sqrt{n}],\ell))\\
&&=\frac{1}{\sqrt{\pi t}}e^{-(x-y)^2/t}-\frac{1-2p_{\ell}}{\sqrt{\pi t}}
e^{-(x+y)^2/t}.
\end{eqnarray*}
\end{thm}

The transition density for Brownian spider in the case of $N=2$, i.e., for
a skew Brownian motion, is given in equations (3) and (4) in Walsh
\cite{WAL} (see also (2.2) in Appuhamillage et al. \cite{AB11}). The transition
density for Brownian spider in the case of $p_j=1/N$,
$j=1,2,\ldots,N,$ is given in Papanicolaou et al. \cite{PPL}.
For general $p_j$, via Walsh, it can be given as follows. Define the
transition density $p(t,v(x,\ell),v(y,\ell^*))$ as
$$
{\bf{P}}({\bb}(t+s)\in v(dy,\ell^*)|{\bb}(s)=v(x,\ell))=
p(t,v(x,\ell),v(y,\ell^*))dy.
$$
As a consequence of Theorem 3.2, we can conclude the following Brownian spider
transition density analogue.
\begin{corollary}
\begin{eqnarray*}
p(t,v(0),v(y,\ell))&=&\frac{2p_\ell}{\sqrt{2\pi t}}e^{-y^2/2t},\\
p(t,v(x,\ell),v(y,\ell^*))&=&\frac{2p_{\ell^*}}{\sqrt{2\pi t}}
e^{-(x+y)^2/2t},\quad \ell\neq \ell^*,\\
p(t,v(x,\ell),v(y,\ell))&=&\frac{1}{\sqrt{2\pi t}}e^{-(x-y)^2/2t}
-\frac{1-2p_\ell}{\sqrt{2\pi t}}e^{-(x+y)^2/2t}.
\end{eqnarray*}
\end{corollary}

\section{Brownian and random walk heights on spider}
\renewcommand{\thesection}{\arabic{section}} \setcounter{equation}{0}
\setcounter{thm}{0} \setcounter{lemma}{0}

One of the natural questions to ask is how high  does the walker go up on
the legs of the spider.
Let $H(j,n)$ denote the highest point reached by the random walk on leg $j$
of the spider in $n$ steps. Formally, let
\begg
\xi(v(r,j),n):=\#\{k: 0<k\leq n,\, \bs_k=v(r,j)\}
\endd
and define
$$H(j,n)=\max\{r:\, \xi(v(r,j),n)\geq 1\}. $$
Let
$$
H_M(n)=\max_{1\leq j\leq N}H(j,n),\qquad H_m(n)=\min_{1\leq j\leq N}H(j,n).$$

Similarly, let ${\bf H}(j,t)$ be the highest point reached by the Brownian
spider ${\bb}(\cdot)$ on leg $j$ by time $t$. Put
$$
{\bf H}_M(t)=\max_{1\leq j\leq N}{\bf H}(j,t),
\qquad {\bf H}_m(t)=\min_{1\leq j\leq N}{\bf H}(j,t).
$$

Note that for fixed $j$ the distribution of $H(j,n)$ and ${\bf H}(j,t)$
can be reduced to the case $N=2$, which is equivalent to skew Brownian
motion and skew random walk. This can be done by keeping the $j-$th leg
as a new leg $1$, and unite all the other legs into leg 2. Then the
distribution of heights $H(j,n)$ and ${\bf H}(j,n)$ are equal to the
distribution of the maximum of skew random walk and maximum of skew
Brownian motion, respectively. The latter one is given in Appuhamillage
and Sheldon \cite{AS12}. Using this result, we obtain the distribution
of ${\bf H}(j,n)$, that also gives  the limiting distribution of $H(j,n)$
as follows.

\begin{thm}
\begg
\lim_{n\to\infty}{\bf P}(H(j,n)<y\sqrt{n})=
{\bf P}({\bf H}(j,t)<y\sqrt{t})=2p_j\sum_{k=1}^\infty (1-2p_j)^{k-1}
(2\Phi((2k-1)y)-1),\label{hdis}
\endd
\textit{where $\Phi$ is the standard normal distribution function.}
\end{thm}

Clearly, $H_M(n)$ and ${\bf H}_M(t)$ are equal to the maximum of a simple
symmetric walk $S(n)$ and of a standard Brownian motion, respectively.
So the law of the iterated logarithm (LIL) and the so called other LIL of
Chung \cite{CH48} continue to hold for these processes.

\begin{thm}
$$
\limsup_{n\to\infty} \frac{H_M(n)}{\sqrt{2n\log\log n}}=
\limsup_{t\to\infty}\frac{{\bf H}_M(t)}{\sqrt{2t\log\log t}}
=1 \quad a.s.
$$
$$
\liminf_{n\to\infty} \left(\frac{\log\log n}{n}\right)^{1/2}H_M(n)
=\liminf_{t\to\infty} \left(\frac{\log\log t}{t}\right)^{1/2}
{\bf H}_M(t)=\frac{\pi}{\sqrt{8}} \quad a.s.
$$
\end{thm}

However, it is a much more interesting question to seek the maximal height
which can be reached on all legs simultaneously. To be more precise, we
are to describe what one can say about $H_m(n)$ and ${\bf H}_m(t)$. For
limsup and Hirsch-type liminf of these processes, we will prove the following
respective results.

\begin{thm}
\begg
\limsup_{n\to\infty} \frac {H_m(n)}{\sqrt{2n\log\log n}}=
\limsup_{t\to\infty} \frac {{\bf H}_m(t)}{\sqrt{2t\log\log t}}=
\frac{1}{2N-1} \quad a.s. \label{twon-1}
\endd

Let $g(t),\, t\geq 1,$ be a nonincreasing function. Then
\begg
\liminf_{n\to\infty}\frac{H_m(n)}{n^{1/2}g(n)}
=\liminf_{t\to\infty}\frac{{\bf H}_m(t)}{t^{1/2}g(t)}=
0\quad or \quad \infty \label{haem}
\endd
according as $\int_1^\infty g(t)\, dt/t$ diverges or converges.
\end{thm}

\noindent
{\bf Proof.}
By the strong approximation given in Section 2, it suffices to prove this
theorem either for $H_m(n)$ or for ${\bf H}_m(t)$.
Denote by
$$
M_1(t)\geq M_2(t)\geq \ldots\geq M_k(t)\geq\ldots
$$
the ranked heights of excursions of a standard Brownian motion on the line
up to time $t$, including the height of a possible incomplete excursion at the
end. It is shown in Cs\'aki and Hu \cite{CsH} it is shown for fixed $k$ that
\begin{equation}
\limsup_{t\to\infty}\frac {M_k(t)}{\sqrt{2t\log\log t}}=
\frac{1}{2k-1} \quad a.s.
\label{mkt}
\end{equation}
and that, for a nonincreasing function $g(t),$
$$
\liminf_{t\to\infty}\frac{M_k(t)}{t^{1/2}g(t)}=
0\quad or \quad \infty
$$
according as $\int_1^\infty g(t)\, dt/t$ diverges or converges.

It is clear that for the Brownian spider, constructed from the standard
Brownian motion as in (\ref{sdef}), we have ${\bf H}_m(t)\leq M_N(t)$,
whenever $M_1(t),\, M_2(t),\, \ldots, M_N(t)$ are on different legs.
Consequently,
$$
\limsup_{t\to\infty}\frac{{\bf H}_m(t)}{\sqrt{2t\log\log t}}
\leq \limsup_{t\to\infty}\frac{M_N(t)}{\sqrt{2t\log\log t}}=
\frac{1}{2N-1}
$$
and
$$
\liminf_{t\to\infty}\frac{{\bf H}_m(t)}{t^{1/2}g(t)}
\leq\liminf_{t\to\infty}\frac{M_N(t)}{t^{1/2}g(t)}=0,
$$
provided $\int_1^\infty g(t)\, dt/t$ diverges.

To show the lower bound in (\ref{twon-1}), let the events $A_n$ and $C_n$ be
defined by
$$
A_n=\left\{M_N(n)\geq (1-\varepsilon)\frac{\sqrt{2n\log\log n}}
{2N-1}\right\},
$$
$$
C_n=\{M_1(n),\, M_2(n),\, \ldots, M_N(n)\,\,
\textrm{are on different legs}\}.
$$
Then, in view of (\ref{mkt}), $P(\limsup_n A_n)=:P(A_n\, \, i.o.)=1$ and,
furthemore, $P(C_n)\geq p_1p_2\ldots p_N=c>0$.

For the next lemma we refer to Klass \cite{KL}.

\begin{lemma}
{\it Let $\{A_n\}_{n\geq 1}$  be an arbitrary sequence of
events such  that  $P(A_n \,\,i.o.)=1$.
 Let $\{C_n\}_{n\geq 1}$ be another arbitrary sequence of events  that is
independent of $\{A_n\}_{n\geq 1}$ and assume that for some $n_0>0$,
$P(C_n)\geq c>0$, for all $n>n_0$.  Then we have}
$P(A_n\,C_n \, \,i.o.)\geq c.$
\end{lemma}

Applying this, we get $P(A_n\,C_n\,\, i.o.)>0$ with $A_n,\, C_n$ as above.
From the $0-1$ law we have also $P(A_n\,C_n\,\, i.o.)=1$. This implies
$$
{\bf P}\left({\bf H}_m(n)\geq (1-\varepsilon)
\frac{\sqrt{2n\log\log n}}{2N-1}\,\,\,\, i.o.\right)=1,
$$
with arbitrary $0<\varepsilon<1$. Hence we have the lower bound in
(\ref{twon-1}).

Now we turn to the convergence part of the liminf result.
From the limit distribution of ${\bf H}(j,t)$ given in (\ref{hdis}), for small
$y$ we have $2\Phi((2k-1)y)-1\leq c(2k-1)y$, hence
$${\bf P}({\bf H}(j,t)<y\sqrt{t})\leq 2cp_jy\sum_{k=1}^\infty
(1-2p_j)^{k-1}(2k-1)=:c_j y,
$$
as $y\to 0$, with some positive constant $c_j$. It is easy to see
that this implies that for any constant $C>0$ we have
\begg
{\bf P}({\bf H}_m(t)<Cy\sqrt{t})<c'y  \label{hirone}
\endd
as well with some positive constant $c'$. By Theorem 2.1, or by the strong
approximation result of Theorem 2.2, we also have
\begin{equation}
{\bf P}(H_m(n)<Cy\sqrt{n})<c'y.
\label{hirone2}
\end{equation}
We prove the convergence part of the liminf in (\ref{haem}) just like that
in Hirsch \cite{H65}. Suppose that $g(n)$ is nonincreasing and
$$
\sum_{n=1}^\infty\frac{g(n)}{n}<\infty.
$$
Then
\begg
\sum_{n=1}^\infty g(2^n)<\infty \label{equi}
\endd
as well. Consequently, from (\ref{hirone2}) we conclude

\begg
{\bf P}\left(\frac {H_m(2^n)}{2^{n/2}}\leq 2C\, g(2^n)\right)\leq 2c'\,g(2^n).
\endd
By (\ref{equi}) and the Borel-Cantelli lemma, we arrive at

\begg
H_m(2^n) \geq 2C\, 2^{n/2}g(2^n) \quad a.s.
\endd
for $n\geq n_0$ with some $n_0$. For an arbitrary $\ell,$ on selecting $k_\ell$
such that
$$2^{k_\ell}<\ell<2\, 2^{k_\ell},$$
 we have
$$H_m(\ell)\geq H_m\left(2^{k_\ell}\right)> 2C\,2^{k_\ell/2}\,
g\left(2^{k_\ell}\right)\geq C\sqrt{\ell}g(\ell)\quad a.s.$$
Since $C$ is arbitrary,
$$
\lim_{n\to\infty}\frac{H_m(n)}{\sqrt{n}g(n)}=\infty\quad a.s.
$$
The convergence part for liminf in (\ref{haem}) is proved. This also completes
the proof of Theorem 4.3. $\Box$

\bigskip
Recall the definitions of $H(j,n)$ and ${\bf H}(j,t)$, i.e., the respective
maximum heights of the random walk on spider and Brownian
spider on leg $j$ up to time $n$ and $t$, respectively.

Our Theorem 4.3. tells us how high could the random walker on a spider, and
Brownian spider, respectively, go up simultaneously on each leg. Now we ask
the following question: if we select N non-negative numbers, as heights,
under what conditions is it possible that the random walker can go up that
high on each leg. The same question can be asked for Brownian spider.
Introducing the notations ${\mathbb R}^N_+$ for the set of vectors with
non-negative components in $N$-dimensional Euclidean space ${\mathbb R}^N$,
i.e.,
$$
{\mathbb R}_+^N:=\{(a(1),\ldots,a(N))\in
{\mathbb R}^N, \, a(1)\geq 0,\ldots,a(N)\geq 0\},
$$
our answer is the following.

\begin{thm}
The set of vectors
\begg
\left(\frac{H(1,n)}{\sqrt{2n\log\log n}},\ldots,
\frac{H(N,n)}{\sqrt{2n\log\log n}}\right),\quad n\geq 3
\label{lim1}
\endd
and
\begg
\left(\frac{{\bf H}(1,t)}{\sqrt{2t\log\log t}},\ldots,
\frac{{\bf H}(N,t)}{\sqrt{2t\log\log t}}\right),\quad t\geq 3
\label{lim2}
\endd
are almost surely relatively compact in ${\mathbb R}^N_+$ and their
respective sets of limit points, as $n\to\infty$ and $t\to\infty$, are given by
\begin{equation}
\left\{(a(1),\ldots,a(N))\in {\mathbb R}^N_+:\, \,
{\cal A}(N):=2\sum_{j=1}^{N}a(j)-\max_{1\leq j\leq N}a(j)\leq 1 \right\}.
\label{condi}
\end{equation}
\end{thm}

For the case $N=2$ equivalent statements are given in Cs\'aki and Hu
\cite{Cshu}, Theorem 1.2, and R\'ev\'esz \cite{R13}, Theorem 5.6.
For the proof we will use the celebrated functional law of the iterated
logarithm of Strassen \cite{ST}. By our strong invariance principle, it suffices
to prove  Theorem 4.4 for random walk on spider. Let ${\cal S}$ be the
Strassen class of functions, i.e., ${\cal S}\subset C([0,1],{\mathbb R})$ is
the class of absolutely continuous functions (with respect to the Lebesgue
measure) on $[0,1]$ for which
\begin{equation}
f(0)=0\qquad {\rm and\qquad } I(f)=\int_0^1\dot{f}^2(x)dx\leq 1.
\end{equation}

Consider the continuous versions of the random walk process
$\{S(nx);\,\, 0\leq x \leq 1\}_{n=1}^{\infty}$ defined  by linear
interpolation from the simple symmetric random walk
$\{S(n)\}_{n=0}^{\infty}.$

\begin{thm} \,\, {\rm \cite{ST}}
The sequence of random functions
$$
\left\{\frac{S(nx)}{(2n\log\log n)^{1/2}};\, 0\leq x\leq
1\right\}_{n\geq 3},
$$
as $n\to\infty$, is almost surely relatively compact in the space
$C([0,1])$ and the set of its limit points is the class of
functions ${\cal S}$.
\end{thm}
{\bf Proof of Theorem 4.4.} Recall the construction of spider walk
from simple symmetric random walk as described in Section 1.
If $a(j)=0$ for all $j=1,2,\ldots,N$, then consider the function $f(x)=0$,
$0\leq x\leq 1$. It is obvious that this function is in $\cal S$, so almost
surely there is a subsequence $n_k$ for which
$$
\lim_{k\to\infty}\sup_{0\leq x\leq 1}\frac{|S(n_kx)|}
{\sqrt{2n_k\log\log n_k}}=0.
$$
This is also true  for the maximums of all excursions. Consequently,
$$
\lim_{k\to\infty}\frac{H(j,n_k)}{\sqrt{2n_k\log\log n_k}}=0,
\quad j=1,2,\ldots,N
$$
i.e., $(0,\ldots,0)$ is almost surely a limit point of (\ref{lim1}) .

Now assume that there are $L$ strictly positive elements among $a(j)$,
$j=1,2,\ldots,N$, denoted by $a(r_1),a(r_2),\ldots, a(r_L)$, and
let $a(r_L)=\max_{1\leq j\leq N}a(j)$ so that we have
\begg
{\cal A}(N)=2\sum_{i=1}^{L-1}a(r_i)+a(r_L)=2\sum_{j=1}^Na(j)
-\max_{1\leq j\leq N}a(j)\leq 1.
\label{boundA}
\endd
We show that $(a(1),\ldots, a(N))$ is almost surely a limit point of
(\ref{lim1}). Construct a piecewise linear function $f(\cdot)$ as follows. Let
$$x_\ell=2(a(r_1)+\ldots+a(r_{\ell-1}))+a(r_\ell),\quad \ell=1,2,\ldots,L$$
and
$$
f(0)=0,\quad f(x_\ell)=(-1)^{\ell-1}a(r_k),\quad \ell=1,2,\ldots,L,\quad
f(1)=f(x_L)=(-1)^{L-1}a(r_L)
$$
and let $f(\cdot)$ be linear in between.

It is easy to see that $f(\cdot)$ is absolutely continuous
and  $I(f)=2\sum_{i=1}^{L-1}a(r_i)+a(r_L)\leq 1,$
consequently $f(\cdot)\in {\cal S}.$ It follows that almost surely there exists
a subsequence $n_k$ such that for the largest $L$ excursion heights
$M(r_i,n_k)$ of $S(n_k)$ we have
$$
\lim_{k\to\infty}\frac{M(r_i,n_k)}{\sqrt{2n_k\log\log n_k}}=a(r_i)\quad a.s.,
$$
for $i=1,2,\ldots,L$ and, if $ M(n_k)$ is another excursion maximum, then we
have
$$
\lim_{k\to\infty}\frac{M(n_k)}{\sqrt{2n_k\log\log n_k}}=0.
$$

For the simple symmetric random walk with $n$ steps, define the event
$$
A_n=\left\{\textrm{there are excursion heights}\, \, M(j,n)\, \,
\textrm{such that}
\, \left|\frac{M(j,n)}{\sqrt{2n\log\log n}}-a(j)\right|\leq
\varepsilon,\, \, j=1,2,\ldots,N\right\}
$$
Then $P(A_n\, i.o.)=1$.

Let $C_n$ be the event that on constructing spider walk from a simple
symmetric random walk, the excursion with height $M(j,n)$ falls to leg $j$
for all $j=1,2,\ldots,N$. $P(C_n)=p_1\ldots p_N>0$, hence, by Lemma 4.1,
$P(A_nC_n \, i.o.)>0$, and by the 0-1 law this probability is 1. Consequently,
$(a(1),a(2),\ldots,a(N))$ is almost surely a limit point of (\ref{lim1}).

To conclude the only if part, assume that $(a(1),\ldots,a(N))$ is a limit point
of (\ref{lim1}), i.e., there exists a subsequence $n_k$ such that
$$
\lim_{k\to\infty}\frac{H(j,n_k)}{\sqrt{2n_k\log\log n_k}}=a(j),\quad
j=1,2,\ldots,N.
$$
Then there exist excursion heights $M(j,n)$ of the random walk $S(i),\,
i=1,2,\ldots,$ and a subsequence $n_k$, $k=1,2,\ldots$ for which
$$
\lim_{k\to\infty}\frac{M(j,n_k)}{\sqrt{2n_k\log\log n_k}}=a(j),\quad
j=1,2,\ldots,N,
$$
and a function $f(\cdot)\in{\cal S}$ such that
$$
|f(x_{2\ell-1})|=a(r_{\ell}),\, \, f(x_{2\ell-2})=0,\, \,
\ell=1,\ldots, L,\, |f(1)|\leq a(r_L),
$$
where, as before, $a(r_1),\ldots, a(r_L)$ are the strictly positive terms
among $a(1),\ldots,a(N)$, and $x_0=0<x_1<x_2<\ldots<x_{2L-1}\leq 1$.

We use the following result (cf. Riesz and Sz.-Nagy \cite{RSzN}, p. 75,
or Shorack and Wellner \cite{SW}, p. 79).
\begin{lemma}
$f(\cdot)\in {\cal S}$ if and only if $f(0)=0$ and for every partition
$0=x_0<x_1<\ldots<x_m=1$ we have
\begg
\sum_{i=1}^m\frac{(f(x_i)-f(x_{i-1}))^2}{x_i-x_{i-1}}\leq 1.
\label{ell2}
\endd
\end{lemma}

This lemma yields
$$
\sum_{\ell=1}^{L-1}a^2(r_\ell)\left(\frac{1}{x_{2\ell-1}-x_{2\ell-2}}
+\frac{1}{x_{2\ell}-x_{2\ell-1}}\right)+\frac{a^2(r_L)}{1-x_{2L-2}}
\label{sum1}
$$
\begg
\leq\sum_{\ell=1}^{L-1}a^2(r_\ell)\left(\frac{1}{x_{2\ell-1}-x_{2\ell-2}}
+\frac{1}{x_{2\ell}-x_{2\ell-1}}\right)+
\frac{a^2(r_L)}{x_{2L-1}-x_{2L-2}} +\frac{(f(1)-a(r_L))^2}{1-x_{2L-1}}
\leq 1.
\label{sum2}
\endd
In case $x_{2L-1}=1$ we take the last term in (\ref{sum2}) to be equal to zero.

The summation of (\ref{sum1}) is of the form
$$
g(z_1,z_2,\ldots,z_{2L-1}):=\sum_{i=1}^{2L-1}\frac{b_i^2}{z_i},
$$
where $z_i>0$ with $\sum_{i=1}^{2L-1}z_i=1$. We want to show that (\ref{sum2})
implies that ${\cal A}(N)\leq 1.$ To this end, we first calculate the minimum
of $g(z_1,z_2,\ldots,z_{2L-1}).$

To find the values of $z_i$ such that the function $g$ takes its minimum, we
have to solve a conditional extreme value problem by the Lagrange multiplier
method, i.e., minimize
$$
g(z_1,\ldots,z_{2L-1})+\lambda(z_1+\ldots+z_{2L-1}-1).
$$
So we have to solve the equations
$$
\frac{b_i^2}{z_i^2}=\lambda,\quad i=1,\ldots,2L-1.
$$
Its solution is $z_i=b_i/(\sum_{i=1}^{2L-1}b_i)$, $i=1,\ldots,L$, i.e., the
minimum value of $g$ is $(\sum_{i=1}^{2L-1}b_i)^2$. Having $g\leq 1,$ by
(\ref{ell2})-(\ref{sum2}) we conclude that $\sum_{i=1}^{2L-1}b_i\leq 1.$
Consequently, for ${\cal A}$ as in (\ref{boundA}), we obtain
$$
{\cal A}(N)=2\sum_{j=1}^{N}a(j)-\max_{1\leq j\leq N}a(j)\leq
2\sum_{i=1}^{L-1}a(r_i)+a(r_L)\leq 1.
$$

This completes the proof of Theorem 4.4. $\Box$

It is worthwhile to give the following corollaries.
\begin{corollary}
Let $M_1(n)\geq M_2(n)\geq\ldots$ be the ranked heights of excursions
of a simple symmetric random walk up to time $n$. Then for finite $N$ we have
$$
\limsup_{n\to\infty}\frac{M_1(n)+2\sum_{i=2}^NM_i(n)}{\sqrt{2n\log\log n}}=1
$$
almost surely.
\end{corollary}
The same is true for Brownian motion.
\begin{corollary}
Let $M_1(t)\geq M_2(t)\geq\ldots$ be the ranked heights of excursions
of a standard Brownian motion up to time $t$. Then for finite $N$ we have
$$
\limsup_{n\to\infty}\frac{M_1(t)+2\sum_{i=2}^NM_i(t)}{\sqrt{2t\log\log t}}=1
$$
almost surely.
\end{corollary}

\section{Increasing number of legs}
\renewcommand{\thesection}{\arabic{section}} \setcounter{equation}{0}
\setcounter{thm}{0} \setcounter{lemma}{0}

In this section we suppose that
$\displaystyle{p_1=p_2=\ldots =p_N=\frac{1}{N}}$.

Let
$$
\xi(v_N(r,j),n):=\#\{k: k\leq n,\, \bs_k=v_N(r,j)\},
$$
i.e., $\xi(v_N(r,j),n)$ is the local time of ${\bf S}$ at time $n$ and
locus $r>0$ on the leg $j$ and, for $r=0$, put
$$
\zeta(n):=\#\{k: k\leq n,\, \bs_k=v_N(0)\}=\xi(v_N(0),n).
$$
Define also the events

\begin{eqnarray*}
M(n,L):&=&\{\min_{1\leq j \leq N} \xi(v_N(L,j),n)\geq 1\} \\
A(n,L,k):&=&\{\min_{1\leq j \leq N} \xi(v_N(L,j),n)\geq k\}.
\end{eqnarray*}
Observe that the meaning of the event $M(n,L)$ is that in $n$ steps the
walker climbs up to at least $L$ on each leg.
The special case $M(n,1)$ means that in $n$ steps each leg is visited at
least once. $A(n,L,k)$ means that in $n$ steps the walker visits each leg
at height $L$ at least $k$ times.

We recall the main result from  R\'ev\'esz \cite{R13}, page 374:

\begin{thm}
For the  {\rm \textbf{SP}}($N$)
\begg
\lim_{N\to \infty} {\bf P}( M([(N\log N)^2],1))=
\left(\frac{2}{\pi}\right)^{1/2} \int_1^\infty e^{-u^2/2}\, du= P(|Z|>1),
\endd
where $Z$  is a standard normal random variable.
\end{thm}

We also have the well-known result that for any $x>0$
\begin{equation}
\lim_{n\to\infty}P\left(\frac{\xi(0,n)}{\sqrt{n}}\leq x\right)=
P(|Z|\leq x).
\label{xin}
\end{equation}

In this section we summerize what we can say about $M(n,L).$ Our main result
is as follows.
\begin{thm} For any integer $L \leq \frac{N}{\log N},$   for
the {\rm \textbf{SP}($N$)} we have
\begg
\lim_{N\to \infty} {\bf P}(M([(cLN\log N)^2],L))=
P\left(|Z|>\frac{1}{c}\right).
\endd
\end{thm}

\smallskip
\noindent
To formulate in words, the theorem above gives the limiting probability of
the event that, as $N\to \infty,$ in $[(cLN\log N)^2]$
steps the walker arrives at least to height $L$ on each of the $N$ legs at
least once.

The next two results are natural companions of the latter one.
\begin{thm}  For any integer $L \leq \frac{N}{\log N},$
{\it and any sequence} $f(N)\uparrow \infty,$   for the
{\rm \textbf{SP}($N$)} we have
\begg
\lim_{N\to \infty} {\bf P}(M([(f(N)L N\log N)^2],L))=1.
\endd
\end{thm}

\begin{thm}  For any integer $L \leq \frac{N}{\log N},$  and any
sequence $f(N)\downarrow 0,$  for the {\rm \textbf{SP}($N$)} we have
\begg
\lim_{N\to \infty} {\bf P}(M([(f(N)L N\log N)^2],L))=0.
\endd
\end{thm}

\noindent
Furthermore we have
\begin{thm} For any integer $L \leq \frac{N}{\log N}$ and any fixed
integer $k\geq 1,$  for the {\rm \textbf{SP}($N$)} we have
\begg
\lim_{N\to \infty} {\bf P}( A([(cLN\log N)^2],L,k))=
P\left(|Z|>\frac{1}{c}\right).
\endd
\end{thm}

\begin{thm} For any integer $L \leq \frac{N}{\log N},$   and any fixed
integer $k \geq 1,$   {\it  and any sequence} $f(N)\uparrow \infty,$
 for the {\rm \textbf{SP}($N$)} we have
\begg
\lim_{N\to \infty} {\bf P}(A([(f(N)L N\log N)^2],L,k))=1.
\endd
\end{thm}

\noindent
{\bf Remark.} In the above five theorems the $L \leq \frac{N}{\log N} $
condition is a technical one, which may be eliminated. So we ask the
following questions.

\noindent
{\bf Question 1:} {\it Determine for each $0\leq p\leq 1,$ the function
$g(N,L,p)$  such that for the {\rm \textbf{SP}($N$)}
$$\lim_{N\to \infty} {\bf P}( M([g(N,L,p)],L))=p$$
should hold.}
\bigskip

\noindent
{\bf Question 2:} {\it Determine for each $0\leq p\leq 1,$ the function
$g^*(N,L,p)$  such that for the {\rm \textbf{SP}($N$)}
$$\lim_{N\to \infty} {\bf P}( A([g^*(N,L,p)],L,k))=p$$
should hold.}

Having assumed in this section that $p_1=p_2=\cdot=p_N=1/N$, in this setup
we can make use of the famous Erd\H{o}s-R\'enyi \cite{ER61} coupon collector
theorem:
\begin{thm} Suppose that there are $N$ urns given, and that
$N\log N+(m-1)N\log\log N+Nx$ balls are placed in these urns one after the
other, independently and equally likely, i.e., with equal probability $1/N$.
Then, for every real $x$, the probability that each urn will contain at
least $m$ balls converges to
\begg
 \exp\left( -\frac{1}{(m-1)!}\exp(-x)\right),
\endd
as
$N\to\infty.$
\end{thm}
It is worthwhile to spell out the most important special case $m=1,$ as
follows.
\begin{thm} Suppose that there are $N$ urns given, and that
$N\log N+Nx$ balls are placed in these urns one after the other,
independently and equally likely, i.e., with probability $1/N$. Then, for
every real $x$, the probability that each urn will contain at least one ball
converges to
\begg
 \exp\left( -\exp(-x)\right),
 \endd
as $N\to\infty.$
\end{thm}

We will also need the following Hoeffding \cite{H61} inequality.
\begin{lemma}
Let $a_i\leq X_i \leq b_i \quad(i=1,2,...k)$
be independent random variables and $S_k=\sum _{i=1}^k X_i.$
Then for every $x>0$
\begg
P(|S_k - E(S_k)|\geq k x)\leq
2\exp\left(-\frac{2k^2x^2}{\sum_{i=1}^k(b_i-a_i)^2}\right). \label{ho1}
\endd
\end{lemma}
We will use the above inequality in the following special case:

Let $X_1, X_2,...X_j\,$ i.i.d. Bernoulli random variables, then for $j\leq k$
\begin{equation}
P(|S_j - E(S_j)|\geq k x)\leq 2\exp\left(-2kx^2\right). \label{ho2}
\end{equation}
To see this, it is enough to observe that for $j\leq k$ we might take
$X_{j+1}=X_{j+2}=...=X_{k}=0$, then $\sum_{i=1}^k(b_i-a_i)^2=j.$

\bigskip
We begin the proofs with some notations.  Let $\{S(n)\}_{n=0}^{\infty}$ be
a  simple symmetric  one-dimensional random walk and let
\begin{eqnarray*}
\xi(0,n)&=&\#\{k: 1\leq k<n,\ S(k)=0\},\\
\zeta(L,n)&:=&\#\{k:1\leq k<n,\ S(k)=0\ {\rm and\,}
|S(k+i)|\,\,\,i=1,2,... \,{\rm  hits}\,\, L \,\,{\rm before \ returning
\ to \,  }\, 0\},\\
\rho(0)&:=&0\, \ {\rm and\,\,} \rho(m):=\min\{k: k>\rho(m-1),\, S(k)=0\},
\end{eqnarray*}
i.e., here, $\zeta(L,n)$ is the number of excursions of the simple symmetric
random walk $S(\cdot)$ reaching $|L|$ before time $n$.

Thus $\xi(0,\rho_m)=m.$ Also observe that $\xi(0,n)=\zeta(1,n)$.

Furthermore, $E(\zeta(L,\rho(i)))=i/L$, on account of
$$
\zeta(L,\rho(i))=\sum_{k=1}^i(\zeta(L,\rho(k))-\zeta(L,\rho(k-1))
$$
being a sum of independent identically distributed Bernoulli random variables
with mean $1/L$. This observation enables us to apply later on the Hoeffding
inequality as above.

Define also
$$ H(n):=\rho(\xi(0,n)+1),$$
i.e., the time of the first return to zero after $n$ steps.

\begin{lemma}
\begin{equation}
|\zeta(L,H(n))-\zeta(L,n)|\leq 1 \quad a.s.
\end{equation}
\begin{equation}
|\xi(0,H(n))-\xi(0,n)|\leq 1\quad a.s.
\end{equation}
\end{lemma}

\vspace{2ex}\noindent
{\bf Proof:} The two statements in hand amounts to observations
in view of the respective definitions of the entities therein.

The next lemma concludes that the probability that the number of excursions
reaching $|L|$ before time $n$ and the number of excursions occuring before $n$
divided by $L$ are too far apart is small.
\begin{lemma}
$$P\left(|\zeta(L,n)-L^{-1}\xi(0,n)|\geq
4 n^{1/4}(\log n)^{3/4}\right)\leq \frac{2}{n}$$
for $n$ big enough.
\end{lemma}
\vspace{2ex}\noindent

\noindent
{\bf Proof:} Let
$$D(n)=|\zeta(L,H(n)) - L^{-1} \xi(0, H(n))|=
              |\zeta(L,\rho(\xi(0,n)+1)) -L^{-1} \xi(0,\rho(\xi(0,n)+1))|.$$
As
$$|\zeta(L,n)-L^{-1}\xi(0,n)|\leq |\zeta(L,H(n))-L^{-1}\xi(0,H(n))|+2,$$
 for $n$ big  enough, we get
\begin{eqnarray*}
\lefteqn{P(|\zeta(L,n)-L^{-1}\xi(0,n)|\geq 4 n^{1/4}(\log n)^{3/4})}\\
&\leq &P(|\zeta(L,H(n))-L^{-1}\xi(0,H(n))|\geq 3 n^{1/4}(\log n)^{3/4})=\\
&= &P(|D(n)|\geq 3\, n^{1/4}(\log n)^{3/4},\,\xi(0,n)\geq2 n^{1/2}
(\log n)^{1/2})+  \\
&+ & P(|D(n)|\geq 3\, n^{1/4}(\log n)^{3/4},\,\xi(0,n)< 2 n^{1/2}
(\log n)^{1/2})=\\
&=& I +II.
\end{eqnarray*}

For $n$ big enough, on account of Lemma 2.2 in Cs\'aki and F\"oldes,
\cite{CsF83} we have
$$
 I\leq P(\xi(0,n)\geq 2 n^{1/2} (\log n)^{1/2})\leq
\frac{C}{n^{2(1-\varepsilon)}}\leq \frac{1}{n},$$
whith an appropriate constant $C>0$ and arbitrary $\varepsilon\in (0,1/2)$.
Furthermore,
\begin{eqnarray*}
II &\leq& \sum_{i=1} ^{2 n^{1/2} (\log n)^{1/2}}
P(|\zeta(L,\rho(i)) -L^{-1}\xi(0,\rho(i))| >
3 n^{1/4} (\log n)^{3/4},\, \xi(0,n)=i )\\
&\leq& \sum_{i=1} ^{2 n^{1/2} (\log n)^{1/2}}
P(|\zeta(L,\rho(i)) -L^{-1}i| >3\, n^{1/4}(\log n)^{3/4})\\
&\leq& \sum_{i=1} ^{2 n^{1/2} (\log n)^{1/2}} 2 \exp(-9\,\log n) \\
& \leq & 4 n^{1/2} (\log n)^{1/2} \exp(-9\,\log n)\leq\exp(-2\log n)=
\frac{1}{n^2}<\frac{1}{n},
\end{eqnarray*}
where we applied Hoeffding inequality  (\ref{ho2}) with
$k=2\, n^{1/2} (\log n)^{1/2}$ and
$\displaystyle{x=\frac{3}{2}\left(\frac{\log n}{n}\right)^{1/4}}.$
$\Box$

\bigskip
\noindent
{\bf Proof of Theorem 5.2} The proof follows the basic ideas of Theorem 5.1.
Suppose that the walker makes $n=[(cLN\log N)^2]$ steps on ${\bf SP}(N).$
This walk can be modelled in the following way. We consider the absolute
value of $S(n),$ where $S(n)$ is a simple symmetric random walk on the line.
Then we get positive excursions which we throw in $N$ urns (the legs of the
spider) with equal probability. We will use Lemma 5.2 to estimate the number
of tall (at least $L$ high) excursions, which are randomly placed in the $N$
urns, and then apply Theorem 5.8. To follow this plan, let

\begin{eqnarray*}
\mu &=& N\log N \\
B_n^{-}&=&\{\zeta(L,n)\leq (1-2\ep) \mu \} \\
B_n&=&\{(1-2\ep)\mu < \zeta(L,n)<(1+2 \ep) \mu \} \\
B_n^+&=&\{  \zeta(L,n)\geq (1+2\ep) \mu \}.
\end{eqnarray*}
In this proof we put $n=[(cLN\log N)^2]$ everywhere, $[\cdot]$ being the
integer part. Having
\begg
\bp(M(n,L))=\bp(M(n,L)|B_n^{-})P(B_n^{-})
+\bp(M(n,L),B_n)+\bp(M(n,L)|B_n^{+})P(B_n^{+}),
\endd
observe that, by Theorem 5.8,
$$\lim_{N\to \infty}\bp(M(n,L)|B_n^{-})=0.$$

\noindent
Using Lemma 5.3,  for $n$ big enough, we have
\begin{eqnarray*}
P(B_n)&=&P((1-2\ep)\mu\leq \zeta(L,n)\leq (1+2\ep)\mu) \\
&\leq&P\left(\frac{\xi(0,n)}{L}\leq(1+3\ep)\mu\right) - \bp
\left(\frac{\xi(0,n)}{L}\leq(1-\ep)\mu\right)
+\frac{4}{n},
\end{eqnarray*}
where we used that the condition $L\log N\leq N$ of the theorem ensures that
\begg
\ep \mu\geq 4 n^{1/4} (\log n)^{3/4} \label{con}
\endd
for large enough $N.$ Consequently, for $N$ big enough, by (\ref{xin}) we have
$$P(B_n)\leq
P\left(|Z|\leq\frac{ (1+3\ep)}{c}\right)-
P\left(|Z|\leq\frac{ (1-\ep)}{c}\right)+o(1).$$
Thus

$$\lim_{N\to \infty}P(B_n)\leq
P\left(|Z|\leq\frac{ (1+3\ep)}{c}\right)-
P\left(|Z|\leq\frac{ (1-\ep)}{c}\right).$$

Again by Theorem 5.8

$$\lim_{N\to \infty}\bp(M(n,L)|B_n^{+})=1,$$
and, by Lemma 5.3, if $N$ is big enough and $L \log N \leq N$,
using (\ref{con}) again, we have that

\begin{eqnarray*}
P(B_n^+)&=&P(  \zeta(L,n)\geq (1+2\ep) \mu ) \\
&\geq&P\left(  \frac{\xi(0,n)}{L}\geq (1+3\ep) \mu \right)+\frac{2}{n}.
\end{eqnarray*}
Consequently, by (\ref{xin}),
$$\lim_{N \to \infty}P(B_n^+)\geq\bp\left(|Z|\geq \frac{1+3\ep}{c}\right)$$
and, similarly,

\begin{eqnarray*}
P(B_n^+)&=&P(\zeta(L,n)\geq (1+2\ep) \mu ) \\
&\leq& P\left(\frac{\xi(0,n)}{L}\geq (1+\ep) \mu \right)+\frac{2}{n}
\end{eqnarray*}
and
$$\lim_{N \to \infty}P(B_n^+)\leq P\left(|Z|\geq \frac{1+\ep}{c}\right).$$
Hence, by the above conclusions, we obtain
$$
P\left(|Z|\geq \frac{1+3\ep}{c}\right)\leq \lim_{N\to\infty}\bp(M(n,L))
$$
$$
\leq P\left(|Z|\leq\frac{ (1+3\ep)}{c}\right)-
P\left(|Z|\leq\frac{ (1-\ep)}{c}\right)+
P\left(|Z|\geq \frac{1+\ep}{c}\right),
$$
for any small enough $\ep>0$.
Letting $\ep\to 0,$ we finally get that
$$\lim_{N \to \infty}\bp(M(n,L))=P\left(|Z|\geq \frac{1}{c}\right).
\quad \quad\Box$$

\bigskip
\noindent
{\bf Proof of Theorem 5.3} We use the  notations of the previous theorem,
with the sole exception that now  $n=[(f(N)L N\log N)^2],$ with
$f(N)\to \infty$.
Observe that
\begg
\bp(M(n,L))\geq\bp(M(n,L)|B_n^{+})\bp(B_n^{+}),
\endd
and, as above, we know that
$$\lim_{N \to \infty}\bp(M(n,L)|B_n^{+})=1.$$
So we only have to show that
$$\lim_{N \to \infty}P(B_n^{+})=1.$$
Now, again by Lemma 5.3,  with any $\ep>0$, for $N$ big enough we have
\begin{eqnarray*}
P(B_n^{+})&=&P(  \zeta(L,n)\geq (1+2\ep) \mu )\geq
P\left(\frac{\xi(0,n)}{L}\geq 4n^{1/4}(\log n)^{3/4}+
(1+2\ep) N \log N\right)- \frac{2}{n} \\
&=&P\left(\frac{\xi(0,n)}{n^{1/2}}\geq
\frac{4n^{1/4}(\log n)^{3/4}}{N\log N f(N)}+
\frac{(1+2\ep) N \log N}{N\log N f(N)}\right)- \frac{2}{n}.
\end{eqnarray*}
Furthermore, having the condition $L\log N\leq N$ and $f(N)\to\infty$, it
is easy to see that
$$\lim_{N \to \infty}\left(\frac{4n^{1/4}(\log n)^{3/4}}{N\log N f(N)}
+\frac{(1+2\ep) N \log N}{N\log N f(N)}\right)=0.$$
Thus the limit of the above probability when $N\to \infty$ is
$P(|Z|\geq 0)=1$, on account of (\ref{xin}). This also proves our
Theorem 5.3. $\Box$

\bigskip
\noindent
{\bf Proof of Theorems 5.5 and 5.6} To prove these two theorems, it is
enough to  repeat the proof of Theorems 5.2 and 5.3, and apply Theorem 5.7
instead of Theorem 5.8.

\bigskip

\noindent
{\bf Proof of Theorem 5.4} Notations are the same as in Theorem 5.3, except
that now $f(N)\downarrow 0$ as $N\to\infty$. During the proof we suppose that
$f(N)L N\log N \to \infty,$  otherwise there is nothing to prove. Observe that
\begg
\bp(M(n,L))\leq\bp(M(n,L)|B_n^{-})P(B_n^{-})+ P(\overline{B_n^{-}}).
\endd
As we know  from Theorem 5.7 that $\lim_{N \to \infty}\bp(M(n,L)|B_n^{-})=0,$
it is enough to prove that
$\lim_{N \to \infty} P(\overline{B_n^{-}})=0.$
 We show that $\lim_{N \to \infty}P(B_n^{-})= 1.$  Using Lemma 5.3 and the
condition $L\log N\leq N,$ with any $0<\ep<1/2$, we have
\begin{eqnarray*} P(B_n^{-})&=&P(\zeta(L,n)\leq (1-2\ep) \mu )
 \\
&\geq&P\left(\frac{\xi(0,n)}{L}+4n^{1/4}(\log n)^{3/4}
\leq (1-2\ep)\mu \right) \\
&\geq&P\left(\frac{\xi(0,n)}{\sqrt{n}}-\frac{2}{n}
\leq \frac{(1-2\ep) }{ f(N)}-
\frac{4n^{1/4}(\log n)^{3/4}}{N\log N f(N)}\right)-\frac{2}{n}\\
&\geq&P\left(\frac{\xi(0,n)}{\sqrt{n}}\leq
\frac{1}{f(N)}\left(1-2\ep-\frac{4f^{1/2}(N) (4 \log N
+2\log f(N))^{3/4}}{\log N}\right)\right)-\frac{2}{n}.
\end{eqnarray*}

Since, now $f(N)\to 0$, as $N\to\infty,$
$\frac{1-2\ep}{f(N)}\to +\infty,$ while the fraction next
to $(1-2\ep)$ goes to 0.
Consequently, by (\ref{xin}), we arrive at
$$\lim_{N\to\infty}P(B_n^{-})=P(|Z|< +\infty)=1.\quad \quad \Box$$

\noindent
{\bf Acknowledgements}

\medskip\noindent
We wish to thank a referee for carefully
reading our manuscript, for posing several probing questions, and providing
also a number of insightful suggestions. All these, in turn, have led to
a much improved presentation of our results. This research was supported by
the Hungarian National Foundation for Scientific Research, Grant No.
K108615, an NSERC Canada Discovery Grant at Carleton University, and by the
PSC CUNY Grant, No. 68030-0043.

\end{document}